\documentclass[11pt]{article}
\usepackage[latin1]{inputenc}
\usepackage[french]{babel}
\usepackage{graphicx}
\usepackage{epsfig}

\usepackage{amsmath}
\usepackage{amssymb}

\newcommand{\titre}[1]{\begin{center}
                      {\Large{\bf #1}}
                      \end{center}}

\newcommand{\orateur}[2]{\begin{center}
{\underline{\large{\bf #1}}}\\
                           {#2}
                          \end{center}}
\begin{document}

\thispagestyle{empty} 

\titre{A derivation of geometric functionals}

\orateur{Thomas Milcent}
{Laboratoire Jean Kuntzmann\\
B.P. 53 - 38041 GRENOBLE Cedex 9 FRANCE \\
e-mail : {\tt thomas.milcent@imag.fr}}


\newcommand{\lra}{\longrightarrow}
\newcommand{\Lra}{\Longrightarrow}
\newcommand{\R}{\mathbb{R}}
\newcommand{\N}{\mathbb{N}}
\newcommand{\Z}{\mathbb{Z}}
\newcommand{\Li}{\mathscr{L}}
\newcommand{\Sp}{\mathbb{S}}
\newcommand{\Ba}{\mathbb{B}}
\newcommand{\Pa}{\mathbb{P}}
\renewcommand{\div}{\operatorname{div}}
\newcommand{\dist}{\operatorname{dist}}
\newcommand{\vect}{\operatorname{vect}}
\newcommand{\Tr}{\operatorname{Tr}}
\newcommand{\Cof}{\operatorname{Cof}}
\newcommand{\rota}{\operatorname{rot}}
\newcommand{\tn}{{\widetilde{n}}}
\newcommand{\tN}{{\widetilde{N}}}
\newcommand{\tC}{{\widetilde{C}}}
\newcommand{\tB}{{\widetilde{B}}}
\newcommand{\tb}{{\widetilde{b}}}
\newcommand{\tc}{{\widetilde{c}}}
\newcommand{\rot}{\nabla\!\times\!}
\newcommand{\Om}{\Omega}
\newcommand{\pOm}{\partial \Omega}
\newcommand{\ds}{\displaystyle}
\newcommand{\ze}{\frac{1}{\eps}\zeta\left(\frac{\phi}{\eps}\right)}
\newcommand{\zep}{\frac{1}{\eps^2}\zeta^{\prime}\left(\frac{\phi}{\eps}\right)}
\newcommand{\zepe}{\frac1{\eps^2}\zeta^{\prime}\left(\frac{\phi}{\eps}\right)}
\newcommand{\gp}{\nabla\phi}
\newcommand{\gd}{\nabla\delta}
\newcommand{\gu}{\nabla u}
\newcommand{\du}{\Delta u}
\newcommand{\gdp}{\nabla \Delta \phi}
\newcommand{\gdu}{\nabla \Delta u}
\newcommand{\ngp}{|\nabla\phi|}
\newcommand{\ngphi}{\frac{\gp}{\ngp}}
\newcommand{\kp}{\kappa(\phi)}
\newcommand{\eps}{\varepsilon}
\newcommand{\phie}{\phi^{\eps}}
\newcommand{\ngpe}{\nabla^{\eps} \phie}
\newcommand{\prpo}{\mathbb{P}_{\gp^\perp}}
\newcommand{\nphi}{\frac{\gp}{\ngp}}
\newcommand{\nsurf}{\nabla_{\pOm}(u\cdot n)}
\newcommand{\nsn}{[\nabla n]\; n}
\newcommand{\tp}{\tilde{\partial}}
\newcommand{\pa}{\partial}
\newcommand{\tsn}{\theta \cdot n}
\newcommand{\usn}{u \cdot n}
\newcommand{\Po}{I\!\!\!P}
\newcommand{\C}{\raisebox{.6ex}
{${\scriptscriptstyle /}$}\hspace{-.43em}C}
\newtheorem{lemme}{Lemma}
\newtheorem{theoreme}{Theorem}
\newtheorem{definition}{Definition}

\section*{Abstract}

We consider in this paper two methods to differentiate functionals of the kind $\int_{\pOm} f(\Om)\; d\sigma$ when $f$ depends only on the geometry of $\pOm$. The first method is based on shape optimization with a Reynolds formula for surfaces. The second consists in a level set approximation on a volumic domain. The main contribution of this paper is to prove that these two approaches lead to the same geometric result when $f$ depends on the normal and the mean curvature and allow to consider general functional depending on the Gaussian curvature.

\section*{Introduction}

An energy functional $J$ depending on a domain $\Om$ is associated in a large range of problems

\begin{equation*}
J(\Om) = \int_{\pOm} f(\Om)\; d\sigma
\end{equation*}

 In order to minimize this energy with a gradient method we have to differentiate $J(\Om)$ with respect to the domain. Shape optimization is now classical. It was introduced originally in applications where $f$ is solution of a PDE \cite{simon,MuratSimon,allaire,Henrot}. We focus here on geometrical shape optimization, in this case  $f$ depends only on the geometry of $\Om$. These energies arise in a lot of applications, for example anisotropic mean curvature flow \cite{} , image processing \cite{} and shape equilibrium of vesicles \cite{}.\\

A  way to perform the geometric shape optimization is to introduce a parametrization of the surface and to compute the derivative in the normal direction. This approach has been used by Willmore  for the mean curvature \cite{Willmore} and by Steigmann  for the Gaussian curvature \cite{Steigmann3} and leads to a geometrical result. However complex tools of differential geometry as the covariant derivative are needed in the computations. To get round of these difficulties we propose two approaches based on a level set representation of the functional.\\

The scope of the paper is the following. In the first section by introducing shape operators and curvature. Then we prove some important lemmas and integration by parts. The second section is devoted to shape optimization of geometrical functionals with a Reynolds formula for surface. The differentiation of functional depending on the normal and the mean curvature are carried out and the result is given in a geometrical form. The third section  deals with another approach based on a volumic approximation with level set functions. Doing so we rewrite the derivative for functionals depending on the normal, mean curvature and Gaussian curvature (with a distance function). We find the results of the second section.

\section{Notations and preliminaries}

In the sequel, $\Om$ will be an regular open connected subset of $\R^3$. Note that its boundary $\pOm$ is a smooth closed surface of $\R^3$. We denote by $\cdot$ the euclidean scalar product on $\R^3$ and $|\cdot|$ the associated norm.  The vectors and matrix will be written in the canonical basis of $\R^3$. For two vectors $a$ and $b$, $a\otimes b$ stands for the tensor product that is $[a\otimes b]_{ij} = a_i b_j$. For two matrix $A,B$ we note $A:B = \sum_{i,j}A_{ij}B_{ij}$ their scalar product. We note $\Tr(A)$, $\Tr(\Cof(A)) = \frac{1}{2}( \Tr(A)^2 - \Tr(A^2))$ and $\det(A)$ the  three invariants  of a matrix $A$ in $M_3(\R)$. We denote by $\nabla$, $\div$ and $\Delta$ the usual gradient, divergence and laplacian, respectively.\\

 We denote by $n$ a unitary smooth extension on $\R^3$ of the normal to the surface $\pOm$. Taking the gradient of $|n|^2=1$ gives

\begin{equation}
[\nabla n]^Tn=0, \hspace{2cm}([\nabla n] \;n)\cdot n = 0.
\label{vp}
\end{equation}

These equations will be used extensively in this paper. Let us first introduce some basics tools of differential geometry as shape operators and curvature. Then we prove some lemmas and formulas to integrate by parts. We conclude this section by introducing the framework of level set methods.

\subsection{Shape operators}

The shape operators generalize classic differential operators on a surface. The reader may consult \cite{Henrot,simon} for a more detailed exposition of this material. The tangential gradient of a smooth function $f : \R^3 \lra \R$ is defined by

\begin{equation*}
\nabla_{\partial \Om} f =(I-n\otimes n)\nabla f=\nabla f- (\nabla f\cdot n)n. 
\label{gradientsurf}
\end{equation*}

The tangential gradient only depends on the values of $f$ on the surface $\pOm$. A simple consequence of the definition is that $\nabla_{\pOm} f \cdot n =0$. The tangential gradient of a smooth vector field $v:\R^3 \lra \R^3$ is the tangential gradient of each component

\begin{equation*}
[\nabla_{\partial \Om} v] = [\nabla v](I-n\otimes n).
\label{gradientsurf2}
\end{equation*}

The tangential divergence of the vector field $v$ is the trace of its tangential gradient

\begin{equation*}
\div_{\pOm}(v) =\Tr([\nabla_{\partial \Om} v]) = \div(v)- ([\nabla v]\; n)\cdot n.
\label{divsurf}
\end{equation*}

The Laplace-Beltrami operator on $\pOm$ of $f$ is the tangential divergence of the tangential gradient

\begin{equation*}
\Delta_{\partial \Om} f = \div_{\partial \Om} (\nabla_{\partial \Om} f).
\label{laplacian}
\end{equation*}
 
These operators are defined on $\R^3$ but they have a geometric interpretation on the surface $\pOm$.

\subsection{Curvature}

The curvature is defined as the variation of the normal in the tangent plane and we denote by $\kappa_1$ and $\kappa_2$ the eigenvalues of the associate endomorphism $dn$. By definition, the mean curvature $H$ is $\kappa_1+\kappa_2$ and the Gaussian curvature $G$ is $\kappa_1 \kappa_2$. The extension of $dn$ on the whole space is the shape operator associate to the normal $[\nabla_{\pOm}n]$. The principal curvatures $\kappa_1$ and $\kappa_2$ are then the nonzero eigenvalues of $[\nabla_{\pOm}n]$. Therefore

\begin{equation}
H= \Tr([\nabla_{\pOm} n]) = \div_{\pOm}(n),\hspace{2cm} G = \Tr(\Cof([\nabla_{\pOm} n])).
\label{courburepOm}
\end{equation}

If the extension of $n$ is unitary we have

\begin{equation}
H= \Tr([\nabla n])  = \div(n), \hspace{2cm} G = \Tr(\Cof([\nabla n])).
\label{courbure}
\end{equation}

The proof is based on (\ref{vp}). Throughout this paper, we assume that the extension of the normal is unitary.

\subsection{Lemmas}

The Leibniz rule to derive the product of functions may be generalized shape operators. Consider two smooth functions $f,g:\R^3 \lra \R$. Following (\ref{courburepOm}) we get

\begin{equation}
\div_{\partial \Om} (f n) = f\div_{\pOm} (n) + \nabla_{\pOm}f\cdot n = f H.
\label{alphan}
\end{equation}

\noindent The operator $P=I-n\otimes n$ is a symmetric projector so

\begin{equation}
\nabla_{\pOm}f\cdot \nabla g=\nabla f\cdot \nabla_{\pOm} g = \nabla_{\partial \Om} f \cdot \nabla_{\pOm} g.
\label{symetrie}
\end{equation}

The normal is unitary so $n=\frac{v}{|v|}$ for some vector $v$. We have the following identity

\begin{equation}
\nsn = \frac{1}{|v|} \left( \nabla |v| - \left(\nabla |v|\cdot \frac{v}{|v|}\right)\frac{v}{|v|}  \right) = \frac{\nabla_{\pOm}|v|}{|v|} 
\label{nsnu}
\end{equation}\\

\noindent We now prove some useful lemmas for the next sections of the paper.

\begin{lemme}
$$\div(\nsn) = \nabla H \cdot n + H^2 - 2G. $$
\label{divnsn}
\end{lemme}

\noindent \textit{Proof}: The property (\ref{courbure}) implies 

\begin{equation*}
\div(\nsn) = \nabla H\cdot n + \Tr([\nabla n]^2)
\end{equation*}

 The lemma is proved with the identity $\Tr([\nabla n]^2) = \kappa_1^2 + \kappa_2^2 = H^2-2G$.\\

\begin{lemme}
$$\left(\nabla \left(\nabla_{\pOm}
f\right) n\right)\cdot n = - \nabla f
 \cdot ([\nabla n]\; n).$$
\label{lemme01}
\end{lemme}

\noindent \textit{Proof}: Taking the gradient of the equality $\nabla_{\pOm}f\cdot n = 0$ and multiplying by $n$ gives

\begin{equation*}
([\nabla(\nabla_{\pOm}f)]^T n) \cdot n + ([\nabla n]^T \nabla_{\pOm} f) \cdot n = 0.
\end{equation*}

\noindent  The lemma is proved combining (\ref{nsnu}) and (\ref{symetrie}).

\begin{lemme}
$$\left(\nabla \left( \nsn \right) n\right)\cdot n = -(\nsn)
 \cdot (\nsn).$$
\label{lemme02}
\end{lemme}

\noindent \textit{Proof}: Taking the gradient of  (\ref{vp}) gives

\begin{equation*}
[\nabla(\nsn)]^T n + [\nabla n]^T \;\nsn=0.
\end{equation*}

\noindent Multiplying by $n$ the result prove the lemma.\\

\subsection{Integration by parts}

Let be $Q$ an open set of $\R^3$. Consider two smooth functions $f,g : \R^3 \lra \R$  and a smooth vector field $v : \R^3 \lra \R^3$. Assume that all these fields vanish on $\partial Q$. We have the following integration by parts for shape operators

\begin{equation}
\int_Q f \div_{\pOm}(v) \; dx=-\int_Q\nabla_{\pOm}f \cdot v\; dx+\int_Q
\div(n\otimes n) \cdot v\;f dx
\label{ipptildev2}
\end{equation}

The proof is a simple consequence of the Stokes formula and the definition of shape operators.\\

We can also integrate by parts shape operators on a closed surface $\pOm$ 

\begin{equation}
\int_{\partial \Om} \nabla_{\partial \Om}f\cdot v =-\int_{\pOm}
f \;\div_{\partial \Om}(v) \;d\sigma +  \int_{\pOm} H f\; v\cdot n \;d\sigma
\label{ipp1}
\end{equation}

The proof of this result can be found in \cite{simon} p 55. We provide another demonstration with a volumic regularization in the third section. Compared to the classical Stokes formula  in $\R^3$ for usual differential operators, an extra term involving the mean curvature of the surface appear. A consequence of (\ref{ipp1}) is the following

\begin{equation}
\int_{\pOm}f\Delta_{\pOm}g \; d\sigma = \int_{\pOm}g\Delta_{\pOm}f\; d\sigma
\label{ipp3}
\end{equation}

\subsection{Level set formulation}

Consider a smooth function $\phi:\R^3\lra \R$ whose $0$ level set represent the smooth surface $\pOm$

\begin{equation*}
\pOm = \{ x\in \R^3 \;/ \;\phi(x)=0  \}
\end{equation*}

Assume that $\ngp>0$ in a neighbourhood of $\{ \phi=0\}$. Consequently, the unitary normal is $n=\frac{\gp}{\ngp}$. If $\phi$ is a distance function ie $\ngp=1$ then $[\nabla n]$ is symmetric and (\ref{vp}) becomes

\begin{equation}
\nsn=0
\label{distance}
\end{equation}\\

In the framework of level set methods, the geometrical functional $\int_{\pOm}f(\Om)\; d\sigma$ may be written

\begin{equation}
J(\phi) = \int_{\{\phi=0\}} f[\phi]\; d\sigma
\label{fonctionnelle}
\end{equation}

where the notation $f[\phi]$ stands for a function  $f$ depending on $\phi$ and its derivatives. A natural idea to perform the shape optimization of $J(\phi)$ is to differentiate with respect to $\phi$. However, we prefer to use an equivalent method which consist in  transporting $\phi$ (which from now on depends on time) with a  velocity field $u(x,t)$

\begin{equation}
\phi_t + u\cdot \gp = 0
\label{transport}
\end{equation}

Therefore, the shape optimization of $J(\phi)$ is performed by differentiating with respect to the parameter $t$. To get round of the problem of differentiating on a domain depending of $\phi$ we propose two methods. The first one consists in applying a Reynolds formula for surfaces and the other one is based on a volumic regularization of the functional.

\section{Optimization with a Reynolds formula for surfaces}

A proof of the Reynolds formula for surfaces is carried out in the annex. Integrating by parts (\ref{reynolds_surface_demo}) with (\ref{ipp1}) gives

\begin{equation}
(J(\phi))_t = \int_{\{\phi=0\}} \; (f[\phi])_t + \div(f[\phi] n)\; u\cdot n\; d\sigma.
\label{result_optiforme}
\end{equation}

 We get the formula obtained in the framework of shape optimization develop by Murat and Simon \cite{MuratSimon}. We focus now on the particular case where $f$ depends on the normal and the mean curvature.

\subsection{Functional depending on the normal}

We first compute the derivative of the normal with (\ref{transport}) and (\ref{nsnu})

\begin{equation}
(n(\phi))_t =\frac{\nabla \phi_t}{\ngp} - \frac{\gp}{\ngp^2} \frac{\nabla \phi_t\cdot \gp}{\ngp}= \frac{-\nabla_{\pOm}  (\usn \;\ngp)}{\ngp} = - \nabla_{\pOm}(u\cdot n) -  (\nsn)\; u\cdot n.
\label{derivee_normale}
\end{equation}

We get the result obtained by shape optimization in \cite{simon} p 70.\\ 

The normal is unitary so it is natural to introduce a smooth function $f:\R^3\lra \R$ homogeneous of degre one and the functional

\begin{equation*}
J^n(\phi) = \int_{\{\phi=0\}} f(n(\phi))\;d\sigma.
\end{equation*}

\noindent By (\ref{result_optiforme}) and the derivative of the normal  (\ref{derivee_normale}) we get

\begin{equation*}
(J^n(\phi))_t = \int_{\{\phi=0\}} \div(f(n) n)\;\usn + \nabla f(n)  \cdot (-\nabla_{\partial \Om}(\usn)- (\nsn)\;\usn)\;d\sigma.
\end{equation*}

\noindent Using the symmetry relation (\ref{symetrie}) and integrating by parts with (\ref{ipp1}) 

\begin{equation*}
- \int_{\{\phi=0\}}\nabla f(n) \cdot \nabla_{\pOm}(\usn) d\sigma  =  -\int_{\{\phi=0\}}\nabla_{\pOm} f(n) \cdot \nabla_{\pOm}(\usn) d\sigma=\int_{\{\phi=0\}} \div_{\pOm}(\nabla_{\pOm} f(n)) \;\usn\; d\sigma.
\end{equation*}

\noindent By (\ref{courbure})

\begin{equation*}
 \div(f(n)n) = f(n) \div(n) + \nabla(f(n))\cdot n = f(n) H + ([\nabla n]^T \nabla f(n)) \cdot n.
\end{equation*}

\noindent Therefore

\begin{equation}
(J^n(\phi))_t = \int_{\{\phi=0\}} (\;f(n)H +\div_{\pOm}(\nabla_{\pOm} f(n)) \;)\; \usn \;d\sigma.
\label{optigene}  
\end{equation}

The relation  $f(tn) = tf(n)$ imply $\nabla f(n)\cdot n = f(n)$. Hence by (\ref{alphan})

\begin{equation*}
\div_{\pOm}( \nabla_{\pOm} f(n)) = \div_{\pOm}(\nabla f(n)) - \div_{\pOm}((\nabla f(n) \cdot n)\; n) =  \div_{\pOm}(\nabla f(n)) - f(n) H.
\end{equation*}

\noindent Finally

\begin{equation}
(J^n(\phi))_t = \int_{\{\phi=0\}} \;\div_{\pOm}(\nabla f(n))\; \usn \; d\sigma.
\label{homogene}
\end{equation}

We get the anisotropic curvature $H_f = \div_{\pOm} (n_f)$ with the anisotropic normal $n_f = \nabla f(n)$ as in \cite{Elliott} p 192.

\subsection{Functional depending on the mean curvature}

\noindent Consider a smooth function $A : \R \lra \R$ and the functional

\begin{equation*}
J^H(\phi)=\int_{\{\phi=0\}}A(H(\phi)) \;d\sigma.
\end{equation*}

The integration by parts formula (\ref{ipp1}) involve shape operators so we choose  the definition (\ref{courburepOm}) of the mean curvature. Following (\ref{result_optiforme})  we get

\begin{equation*}
(J^H(\phi))_t = \int_{\{\phi=0\}} \div(A(H) n)\;\usn + A^{\prime}(H) (H(\phi))_t \;d\sigma.
\end{equation*}

\noindent The derivative of the normal (\ref{derivee_normale}) gives

\begin{eqnarray*}
(H(\phi))_t =&& 
\div(-\nsurf -(\nsn)\; \usn) - (\nabla (-\nsurf -(\nsn)\; \usn)\; n) \cdot n\\
&&-  ([\nabla n]\; (-\nsurf -(\nsn) \;\usn) )\cdot n-  (\nsn) \cdot (-\nsurf -(\nsn) \;\usn).
\end{eqnarray*}

\noindent We put the two first terms in a tangential divergence and the third one vanish by (\ref{vp}). Developping the divergence gives

\begin{equation*}
(H(\phi))_t= - \Delta_{\pOm}(\usn) - \div_{\pOm}(\nsn)\; \usn   + (\nsn) \cdot(\nsn) \;\usn.
\end{equation*}

\noindent  Following  lemma \ref{divnsn} and lemma \ref{lemme02}

\begin{equation*}
\div_{\pOm}(\nsn) = \div(\nsn) - (\nabla(\nsn) n)\cdot n =  \nabla H \cdot n + H^2 -2G + (\nsn)\cdot (\nsn).
\end{equation*}

\noindent We have also $\div(A(H)n) = A(H) H + A^{\prime}(H)\nabla H\cdot n$. After integrating by parts the laplacian with (\ref{ipp3}) we finally get

\begin{equation}
(J^H(\phi))_t = \int_{\{\phi=0\}} (A(H)H - \Delta_{\pOm}(A^{\prime}(H)) 
 -A^{\prime}(H)(H^2 -2G)) \; \usn \; d\sigma.
\label{resopti}
\end{equation}

We get the result of Willmore \cite{Willmore} which is geometrical because it involves only the curvature and shape operators. For a more general functional depending on the Gaussian curvature, complex computations gives the same result as the one developed in next section devoted to optimization by volumic approximation.

\section{Optimization with a volumic approximation}

The following lemma is the key of method of optimization by volumic approximation.

\begin{lemme}
Let be $\zeta$ a continuous positive cut-off function with support in $[-1,1]$ and $\int_{\R}\zeta=1$. Assume that $f,\phi$ are regular and $\ngp>0$, then
$$J_{\eps}(\phi) := \int_Q f[\phi]\ngp \ze\; dx \underset{\eps \lra 0}{\lra} \int_{\{ \phi=0\}} f[\phi]\; d\sigma$$
\label{approx}
\end{lemme}

\noindent \textit{Proof}: The  co-area formula is

\begin{equation}
\int_{\{|\phi(x)|<\eps\}} f(x) \; dx = \int_{-\eps}^\eps \int_{\{\phi(x)=\nu\}} f(x)|\gp(x)|^{-1}\;d\sigma d\nu
\label{coaera}
\end{equation}

This formula is quite natural since the volume $\{ |\phi|<\eps\}$ is computed integrating the area of $\{ \phi=\nu\}$ with a factor $\ngp^{-1}$ which correspond to the spacing between the level set of $\phi$. We refer to \cite{Evans} p 118 for a rigorous proof. For a continuous function $g$ and under the hypothesis on $\zeta$ we have $\int_{-\eps}^{\eps} \frac{1}{\eps} \zeta\left(\frac{r}{\eps}  \right) g(r)\;dr \underset{\eps\lra 0}{\lra} g(0)$. Applying the co-area formula (\ref{coaera}) to the continuous function $g(r) = \int_{\{\phi(x)=r\}} f(x) \; d\sigma$ prove the lemma.\\

{\bf Remark}

This lemma allows us to get the integration by parts (\ref{ipp1}). Indeed, the integration by parts for shape operators (\ref{ipptildev2}), the property $\nabla_{\pOm}\phi=0$ and (\ref{nsnu}) give

\begin{equation*}
\int_Q \div_{\pOm}(v) \;\ngp \ze\; dx =
\int_Q H \;v\cdot n\;\ngp \ze\; dx.
\end{equation*}

\noindent Taking $\eps\lra 0$ with the lemma \ref{approx} and replace $v$ by $fv$ gives the result.

\subsection{General functionnal}

We compute the derivative of $J_{\eps}(\phi)$ with respect to $t$

\begin{equation*}	
(J_{\eps}(\phi))_t =\int_Q (f[\phi])_t \;\ngp \ze \; dx  + \int_Q f[\phi]\ze
\frac{\gp}{\ngp} \cdot \nabla \phi_t \;dx + \int_Q f[\phi] \ngp \frac{1}{\eps^2}\zeta^{\prime}\left(\frac{\phi}{\eps}\right) \;\phi_t\;dx.
\end{equation*}

\noindent Integrating by parts the second term gives

\begin{equation*}
\int_Q f[\phi] \ze
\frac{\gp}{\ngp} \cdot \nabla \phi_t\; dx= -\int_Q
\div\left(f[\phi]\frac{\gp}{\ngp}\right)\ze \phi_t\; dx -\int_Q
f[\phi]\frac{\gp}{\ngp} \cdot \gp\;\frac{1}{\eps^2}\zeta^{\prime}\left(\frac{\phi}{\eps}\right) \phi_t\ \;dx.
\end{equation*}	

\noindent Hence (\ref{transport}) implies

\begin{equation}
(J_{\eps}(\phi))_t =  \int_Q \left(\;(f[\phi])_t +  
\div\left(f[\phi]n\right) u\cdot n \;\right) \ngp \ze  \;dx.
\label{surflevel1}
\end{equation}

We get the Reynolds formula for surfaces (\ref{result_optiforme}) taking $\eps\lra 0$ with the lemma \ref{approx}. We focus now on the case where $f$ depends on the normal and the curvature.

\subsection{Functional depending of the normal}

Consider a smooth function $f:\R^3 \lra \R$ homegenous of  degree one and the functional

\begin{equation*}
J_{\eps}^n(\phi) = \int_{Q} f\left(n(\phi)\right)\ngp \ze\;dx.
\end{equation*}

 Following (\ref{surflevel1}) and the derivative of the normal (\ref{derivee_normale}) we obtain

\begin{equation*}
(J_{\eps}^n(\phi))_t = \int_{Q} \div(f(n) n)\;\usn \;\ngp \ze \; dx  + 
\int_{Q}\nabla f(n)\cdot \nabla_{\pOm}\phi_t\ze\;dx.
\end{equation*}

\noindent  Using the symetric relation (\ref{symetrie}) and  integrating by parts gives

\begin{equation*}
 \int_{Q}\nabla f(n)\cdot \nabla_{\pOm}\phi_t\ze\;dx =   \int_{Q}\nabla_{\pOm} f(n) \cdot \nabla \phi_t\ze\;dx= - 
\int_{Q}\div(\nabla_{\pOm} f(n))\ze \phi_t\;dx.
\end{equation*}

\noindent (\ref{transport}) implies

\begin{equation}
(J_{\eps}^n(\phi))_t = \int_{Q} \div(\;f(n) n  +
\nabla_{\pOm} f(n)\;) \;\usn \;\ngp \ze \;dx.
\label{levelgene}
\end{equation}

\noindent This formula may be written in another form. Indeed with lemma \ref{lemme01}

\begin{equation*}
 \div(\nabla_{\pOm} f(n)) = \div_{\pOm}(\nabla_{\pOm} f(n)) + (\nabla(\nabla_{\pOm} f(n))\; n)\cdot n= \div_{\pOm}(\nabla_{\pOm} f(n)) -\nabla f(n) \cdot (\nsn).
\end{equation*}

\noindent We have also $\div(f(n) n ) = f(n) H  + ([\nabla n]^T\nabla f(n))\cdot n$. Therefore we get the previous result (\ref{optigene}) taking $\eps\lra 0$ with lemma \ref{approx}.

\subsection{Functional depending of the curvature}

Consider a smooth function $g:M_3(\R)\lra \R$ and the general functionnal

\begin{equation*}
J_{\eps}(\phi) = \int_Q g\left(\nabla (n(\phi)) \right)\ngp \ze \; dx.
\end{equation*}

We denote by $[\nabla g]$ the gradient of $g$ at the point $[\nabla n]$. Following (\ref{surflevel1}) and the derivative of the normal (\ref{derivee_normale})

\begin{equation*}
(J_{\eps}(\phi))_t = \int_{Q} \div(g([\nabla n]) n )\;\usn\;\ngp\ze \;dx + \int_{Q} [\nabla
g]:\nabla\left(\frac{\nabla_{\pOm} \phi_t}{\ngp}\right) \ngp \ze \; dx.
\end{equation*}

\noindent Integrating by parts the second term gives

\begin{equation*}
\int_{Q} [\nabla
g]:\nabla\left(\frac{\nabla_{\pOm} \phi_t}{\ngp}\right) \ngp \ze \; dx
= -  \int_{Q}\div\left( [\nabla g]\ngp \ze\right)\cdot
 \frac{\nabla_{\pOm} \phi_t}{\ngp} \; dx.
\end{equation*}

\noindent The integration by parts (\ref{ipptildev2}) for shape operators gives

\begin{eqnarray*}
(J_{\eps}(\phi))_t
= &&\int_{Q} \div(g([\nabla n]) n)\; \usn\;\ngp\ze \;dx+ \int_{Q} \div\left( [\nabla g]\ngp\ze\right)\cdot \div(n\otimes n)\;\usn\;dx\\ &-& \int_{Q}\div_{\pOm}\left(\frac{1}{\ngp} \div\left( [\nabla g]\ngp\ze\right)\right)
 \;\usn \; \ngp\; dx
\end{eqnarray*}

Developping the divergence

\begin{equation*}
\div\left( [\nabla g]\ngp\ze\right) = \div( [\nabla g]) \ngp \ze + [\nabla g]\;\nabla\ngp \ze + [\nabla g]\;\gp \;\ngp \zepe.
\end{equation*}

\noindent The relation $\nabla_{\pOm}(\phi) = 0$ implies $\div_{\pOm}\left(\ze v\right)= \div_{\pOm}(v)\ze$ for all vector $v$. We obtain finally

{\small
\begin{eqnarray*}
(J_{\eps}(\phi))_t &=&  \int_{Q} \div(g([\nabla n]) n )\;\usn\;\ngp \ze \;dx
 - \int_Q \div_{\pOm}(\div( [\nabla g]))\;\usn\;\ngp  \ze \;  \;dx\\
 &+& \int_Q \div( [\nabla g])\cdot\div(n\otimes n)  \;\usn\;\ngp \ze  \;dx
+ \int_Q\left([\nabla g]\; \frac{\nabla\ngp}{\ngp}\right)\cdot\div(n\otimes n)\;\usn\;\ngp \ze \;dx\\
&-&\int_Q\div_{\pOm}\left([\nabla g]\;\frac{\nabla\ngp}{\ngp}\right)\;\usn\;\ngp \ze\;  \;dx
 -\int_Q\div_{\pOm}\left([\nabla g]\; \gp \zepe\right) \;\usn\;\ngp \;dx\\
&+&\int_Q\left([\nabla g]\; \gp \right)\cdot\div(n\otimes n)\zepe\; \;\usn\;\ngp  \;dx.
 \label{levfin}
\end{eqnarray*}
}

In order the simplify this complex expression in the case of a functional depending on the mean and Gaussian curvature we consider 

\begin{equation*}
 g(w) = F(\Tr(w), \Tr(\Cof(w) )).
\end{equation*}

where $F:\R\times \R \lra \R$ is a smooth function. We denote $F_H$ and $F_G$ the derivative of $F$ with respect to his first and second variable, respectively. We have

\begin{equation*}
 [\nabla g]([\nabla n]) =  \alpha I + \beta[\nabla n]^T.
\end{equation*}

where $\alpha = F_H + H F_G$ and $\beta = -F_G$.\\

 The two last terms involving $\zeta^{\prime}$ in the expression of $(J_{\eps}(\phi))_t$ seem annoying. However using (\ref{vp}) and (\ref{alphan}) these terms compensate.\\

\noindent We now simplify the fourth and fifth terms involving $\frac{\nabla\ngp}{\ngp}$. Following (\ref{nsnu}) we have $\frac{\nabla\ngp}{\ngp} = \nsn + \frac{\nabla\ngp\cdot n}{\ngp} n$. Hence with (\ref{alphan}) and (\ref{vp})

\begin{equation*}
-\div_{\pOm}\left((\alpha I + \beta [\nabla n]^T) \frac{\nabla\ngp}{\ngp}\right) 
=-\div_{\pOm}\left((\alpha I + \beta [\nabla n]^T)\nsn\right) 
-\alpha H \frac{\nabla\ngp\cdot n}{\ngp}.
\end{equation*}

\noindent We have also

{\small
\begin{equation*}
 \left((\alpha I + \beta [\nabla n]^T)\; \frac{\nabla\ngp}{\ngp}\right)\cdot\div(n\otimes n) =(\alpha + \beta H)(\nsn)\cdot (\nsn) +\alpha H\frac{\nabla\ngp \cdot n}{\ngp} + \beta ([\nabla n]\;
\nsn)\cdot(\nsn).
\end{equation*}
}

\noindent We denote by $I_{\eps}$ the expression of $(J_{\eps}(\phi))_t$ without the integrals and the terms $\usn \ngp \ze$. Therefore

\begin{eqnarray*}
I_{\eps}= &&\div(g([\nabla n]) n)- \div_{\pOm}(\div( \alpha I + \beta [\nabla n]^T)) 
 - \div_{\pOm}\left((\alpha I + \beta [\nabla n]^T)\; \nsn\right)\\
&+& \div( (\alpha I + \beta [\nabla n]^T))\cdot\div(n\otimes n) 
 +(\alpha + \beta H)(\nsn)\cdot (\nsn)+ \beta ([\nabla n]\;
\nsn)\cdot(\nsn).
\label{levs2} 
\end{eqnarray*}

We propose to compute the contributions of $\alpha I$ and $\beta [\nabla n]^T$ in order to obtain the derivative of a functional depending of the mean and the Gaussian curvature.

\subsection*{Functional depending on the mean curvature}

We note $R^{\alpha}$ the expression $I_{\eps}$ with $\beta=0$ and without the term $\div(g([\nabla n]) n)$

\begin{equation*}
R^{\alpha}= - \div_{\pOm}(\nabla \alpha) 
 - \div_{\pOm}\left(\alpha\nsn\right)
+ \nabla \alpha\cdot(H n + \nsn)
 +\alpha (\nsn)\cdot (\nsn).
\end{equation*}

\noindent (\ref{alphan}) implies $\div_{\pOm}(\nabla \alpha) = \Delta_{\pOm} \alpha + H\nabla \alpha \cdot n$. Following lemma \ref{lemme02}, (\ref{symetrie}) and (\ref{nsnu})

\begin{equation*}
- \div_{\pOm}(\alpha\nsn) = -\alpha \div(\nsn)- \alpha (\nsn) \cdot (\nsn)- \nabla \alpha\cdot(\nsn).
\end{equation*}

\noindent We finally get

\begin{equation}
R^{\alpha}= - \Delta_{\pOm}\alpha - \alpha \div(\nsn).
\label{levalpha2}
\end{equation}

Consider the functional 

\begin{equation*}
J_{\eps}^H(\phi)=\int_{Q}A(H(\phi))\ngp \ze\; dx
\end{equation*}

  which correspond to $\alpha=A^{\prime}(H)$ and $\beta=0$.  The relation $\div(A(H)n)=A(H)H + A^{\prime}(H)\nabla H\cdot n$ and the lemma \ref{divnsn} gives

\begin{equation*}
(J^H_{\eps}(\phi))_t =\int_{Q}\left( \; A(H)H - \Delta_{\pOm}(A^{\prime}(H)) - A^{\prime}(H) (H^2-2G)\;\right) \;\usn\;\ngp \ze dx.
\end{equation*}

We get the result (\ref{resopti}) taking $\eps\lra 0$ with the lemma \ref{approx}.

\subsection*{Functional depending on the Gaussian curvature with a distance function}

 We choose a distance function for $\phi$ to simplify $I_{\eps}$. We note $R^{\beta}_d$ the expression of $I_{\eps}$ with  $\alpha =0$  and without the term $\div(g([\nabla n])n)$ in the case of a distance function. Therefore, by (\ref{distance})

\begin{equation*}
R^{\beta}_d= - \div_{\pOm}(\div( \beta [\nabla n]^T)) 
+ \div( ( \beta [\nabla n]^T))\cdot H n.
\end{equation*}

\noindent Following (\ref{distance}) we get $\div(\beta [\nabla n]^T)\cdot H n = \beta H\nabla H\cdot  n$. We have also with (\ref{vp}), (\ref{alphan}) and the property of symmetry (\ref{symetrie})

\begin{equation*}
-\div_{\pOm}(\div( \beta [\nabla n]^T)) = -\beta\Delta_{\pOm} H - \beta H (\nabla H \cdot n)- \nabla_{\pOm} H \cdot \nabla_{\pOm}\beta - \div_{\pOm}([\nabla n]^T\; \nabla_{\pOm}\beta).
\end{equation*}

\noindent The relation $\div_{\pOm}(Ab) = \div_{\pOm}(A^T)\cdot b + A^T : \nabla_{\pOm} b$ and the symmetry of  $[\nabla n]$ gives

\begin{equation*}
-\div_{\pOm}(  [\nabla n]^T\;\nabla_{\pOm}\beta ) = -\div_{\pOm}(  [\nabla n]^T)\cdot\nabla_{\pOm}\beta  - [\nabla n] : \nabla_{\pOm}(\nabla_{\pOm}\beta ).
\end{equation*}

\noindent Taking the gradient of (\ref{vp}) and multiply by $n$ gives $(\nabla([\nabla n]^T) n) \cdot n = -[\nabla n]^T \nsn = 0$. Therefore $\div_{\pOm}(  [\nabla n]^T) = \nabla H$. The symetric relation (\ref{symetrie}) and the relation $[\nabla n] = [\nabla_{\pOm} n]$ valid for a distance function gives

\begin{equation}
R^{\beta}_d=- \beta \Delta_{\pOm} H -2 \nabla_{\pOm} H\cdot \nabla_{\pOm}\beta - [\nabla_{\pOm} n] : \nabla_{\pOm}(\nabla_{\pOm}\beta ).
\label{Rbetad}
\end{equation}\\

Consider the general functional

\begin{equation*}
J_{\eps}^{H,G}(\phi) = \int_Q F(H(\phi),G(\phi))\ngp \ze dx
\end{equation*}

 Let $J_d$ the sum of $R^{\alpha}$ (\ref{levalpha2}), $R^{\beta}_d$ (\ref{Rbetad}) and the term $\div(F(H,G)n)$ with $\alpha = F_H + H F_G$ and $\beta = - F_G$. In case of a distance function we have $R^{\alpha}= -\Delta_{\pOm}\alpha$ so

\begin{equation*}
J_d= \div(F(H,G)n)- \Delta_{\pOm}(F_H + H F_G) + F_G \Delta_{\pOm}H + 2\nabla_{\pOm}H\cdot \nabla_{\pOm}F_G
+  [\nabla_{\pOm} n]:\nabla_{\pOm}( \nabla_{\pOm}F_G  ). 
\end{equation*}

\noindent Developing the surfacic laplacian

\begin{equation*}
-\Delta_{\pOm}(H F_G) + F_G \Delta_{\pOm}H + 2\nabla_{\pOm}H\cdot \nabla_{\pOm}F_G= -H \Delta_{\pOm}F_G.
\end{equation*}

\noindent Using the lemma \ref{divnsn} with a distance function

\begin{equation*}
 \div(F(H,G)n) = F(H,G)H - F_H\; (H^2-2G) +  F_G\; \nabla G \cdot n.
\end{equation*}

\noindent Taking $\eps\lra 0$ with the lemma \ref{approx} gives

\begin{eqnarray}
(J^{H,G}(\phi))_t = &\;& \int_{\{\phi=0\}}\left(\;F(H,G)H - F_H\;(H^2 - 2G) -  \Delta_{\pOm}F_H \;\right)\usn\; d\sigma \\
 &+& \int_{\{\phi=0\}}\left(\;F_G\; \nabla G \cdot n - H\Delta_{\pOm}F_G 
+ [\nabla_{\pOm} n]:\nabla_{\pOm}( \nabla_{\pOm} F_G )\;\right) \usn\; d\sigma
\label{steigmann2}
\end{eqnarray}

\noindent Assume the surface do not change of topology. The Gauss-Bonnet theorem claims that $\int_{\{\phi=0\}} G\; d\sigma$ is constant. Applying the formula (\ref{steigmann2}) with $F(H,G)=G$ gives ($F_H=0$ and $F_G=1$) \; $\int_{\{\phi=0\}} (GH + \nabla G \cdot n) \;\usn\; d\sigma = 0$. This relation is valid for all $u$ so $\nabla G \cdot n = -GH$ on $\{\phi=0\}$

\noindent We obtain finally

\begin{eqnarray}
(J^{H,G}(\phi))_t = &\;& \int_{\{\phi=0\}}\left(\;F(H,G)H - F_H\;(H^2 - 2G) -  \Delta_{\pOm}F_H \;\right)\usn\; d\sigma \\
 &+& \int_{\{\phi=0\}}\left(\;-F_G\;GH - H\Delta_{\pOm}F_G 
+ [\nabla_{\pOm} n]:\nabla_{\pOm}( \nabla_{\pOm}F_G )\;\right) \usn\; d\sigma
\label{steigmann3}
\end{eqnarray}

We get the result of Steigmann \cite{Steigmann3} which is geometric because it involve only shape operators and curvature.

\section{Conclusion}

In this paper two methods of shape optimization had been investigated for functionals depending on the normal and the curvature. The first one is based on a Reynolds formula for surfaces and integration by parts forces us to introduce shape operators. The second one, based on volumic approximation with level set leads to simple computations but some efforts are needed to put the results on a geometrical form. These two approaches give the same result which is geometrical for functional depending on the normal and the mean curvature. The second methods allow us to consider more general functionals depending on the Gaussian curvature. These formulas will be useful to study the influence of the Gaussian curvature in the shape equilibrium of vesicles.

\section{Annex}

To be self contained we present here a demonstration of the Reynolds formula for surfaces based on change of variables on surfaces.\\

\subsection{Change of variables in surface integral}

\noindent Let be $U$ an open set of $\R^2$ and a smooth function $\Phi:\R^3\times \R\lra \R$.  Consider $X: U \lra \R^3$ a smooth  parametrisation of the surface $\pOm_0$ and $Y:U\times \R \lra \R^3$ a parametrisation associate to the surface $\pOm_t=\Phi(\pOm_0,t)$. Consider also a smooth function $f:\R^3\times \R\lra \R$. We have

\begin{equation*}
\int_{\pOm_t} f(x,t)\; d\sigma = \int_U f(Y(u^1,u^2,t)) |Y_{,1} \wedge Y_{,2}|\;du^1 du^2.
\end{equation*}

The relation $Y_{,i} = [\nabla \Phi] \;X_{,i}$ gives for all vector $v$

\begin{equation*}
(Y_{,1} \wedge Y_{,2})\cdot ([\nabla \Phi] \;v)  = \det( [\nabla \Phi]\; X_{,1},[\nabla \Phi]\; X_{,2},[\nabla \Phi]\; v) = \det([\nabla \Phi]) \det(X_{,1},X_{,2},v).
\end{equation*}

Therefore

\begin{equation}
Y_{,1} \wedge Y_{,2}= \det([\nabla \Phi]) [\nabla \Phi]^{-T}X_{,1}\wedge X_{,2}.
\label{norm_inter}
\end{equation}
 
The relations $n_0 = \frac{X_{,1}\wedge X_{,2}}{|X_{,1}\wedge X_{,2}|}$ and $n = \frac{Y_{,1}\wedge Y_{,2}}{|Y_{,1}\wedge Y_{,2}|}$ implies
$|Y_{,1} \wedge Y_{,2}|  = |\Cof([\nabla \Phi]) n_0| |X_{,1}\wedge X_{,2}|$  so

\begin{equation}
\int_{\pOm_t} f(x,t)\; d\sigma =  \int_{\pOm_0} f(\Phi(a,t),t)|\Cof([\nabla \Phi](a,t)) n_0(a)|\; d\sigma_0.
\label{ch_variable_surface3}
\end{equation}

We have also from (\ref{norm_inter})

\begin{equation}
n = \frac{\Cof([\nabla \Phi]) n_0}{|\Cof([\nabla \Phi]) n_0|}. 
\label{normale_phiat}
\end{equation}

\subsection{Reynolds formula for surfaces}

Consider a surface $\pOm_t$ evolving in a smooth velocity field $u:\R^3\times \R\lra \R^3$ and $f:\R^3\times \R\lra \R$ a smooth  function. Then we have the following Reynolds formula for surfaces

\begin{equation}
\frac{d}{dt}\left(\int_{\pa \Omega_t} f(x,t) \; d\sigma\right) = \int_{\pa \Omega_t} f_t + u\cdot \nabla f + f \div_{\pOm}(u)\; d\sigma.
\label{reynolds_surface_demo}
\end{equation}

\noindent \textit{Proof}: The main idea consist to write the integral on a fixed domain $\pOm_0$ and differentiate with respect to $t$. Consider the caracteristics $\Phi$ associates to the velocity field $u$ ie $\frac{\pa \Phi}{\pa t}(a,t) = u(\Phi(a,t),t)$ and denote $\pOm_t = \Phi(\pOm_0,t)$. The change of variables $x= \Phi(a,t)$ (\ref{ch_variable_surface3}) gives

\begin{equation*}
\int_{\pOm_t} f(x,t)\; d\sigma = \int_{\pOm_0} f(\Phi(a,t),t) \det([\nabla \Phi(a,t)])|[\nabla \Phi(a,t)]^{-T} n_0(a)|\; d\sigma_0.
\end{equation*}

Differentiating with respect to t gives three terms. The caracteristics equation gives for the first term $(f(\Phi,t))_t = f_t + u\cdot \nabla f$. The caracteristics equation gives also 

\begin{equation}
 [\nabla \Phi]_t = [\nabla u][\nabla \Phi]
\label{nablaphit}
\end{equation}

The derivative of the determinant $(\det(A(t)))^{\prime} = \det(A(t))\Tr([A(t)]^{-1}A^{\prime}(t))$ and (\ref{nablaphit}) gives for the second term 

\begin{equation*}
(\det(\nabla \Phi))_t = \det(\nabla \Phi) \div(u).
\end{equation*}

\noindent  For the third term, the derivative of the inverse  $(A^{-1})^{\prime}(t) = -A^{-1}(t) A^{\prime}(t)A^{-1}(t)$, (\ref{nablaphit}) and (\ref{normale_phiat}) gives

\begin{equation*}
|[\nabla \Phi]^{-T} n_0|_t = (([\nabla \Phi]^{-T})_t n_0)\cdot \frac{[\nabla \Phi]^{-T} n_0}{|[\nabla \Phi]^{-T} n_0|}= -\left( [\nabla u]\; n \right)\cdot n \;|[\nabla \Phi]^{-T} n_0|.
\end{equation*}

\noindent The change of variables $x=\Phi(a,t)$ and  the definition of tangential divergence prove the lemma.

\bibliographystyle{plain}
\bibliography{bibliographie_thomas}

\end{document}